\documentclass[12pt]{article}
\usepackage{amsmath}
\usepackage{amssymb}
\usepackage{amsthm}
\usepackage{cite}
\usepackage{enumerate}
\def\R{\mathbb{R}}
\def\Z{\mathbf{Z}}
\def\N{\mathbb{N}}

\newcommand{\nc}{\newcommand}
\nc{\BC}{\mathbb{C}}
\nc{\BS}{\mathbb{S}}
\nc{\BP}{\mathbb{P}}
\nc{\BE}{\mathbb{E}}
\nc{\BQ}{\mathbb{Q}}
\nc{\bN}{{\mathbf N}}
\nc{\BX}{{\mathbb X}}
\nc{\bX}{{\mathbf X}}
\nc{\bY}{{\mathbf Y}}
\nc{\BY}{{\mathbb Y}}
\nc{\bM}{{\mathbf M}}
\nc{\bF}{{\mathbf F}}
\nc{\bG}{{\mathbf G}}
\nc{\bH}{{\mathbf H}}
\nc{\bU}{{\mathbf U}}
\nc{\bz}{{\mathbf z}}
\nc{\bx}{{\mathbf x}}
\nc{\be}{{\mathbf e}}
\nc{\bW}{{\mathbf W}}
\DeclareMathOperator{\BV}{{\mathbb Var}}
\DeclareMathOperator{\CV}{{\mathbb Cov}}

\DeclareMathOperator{\EV}{{\mathbb E}}
\DeclareMathOperator{\cl}{{\rm cl}}

\begin{document}

\newcommand{\I}{{\bf 1}}
\newcommand{\CF}{{\cal F}}
\newcommand{\CA}{{\cal A}}
\newcommand{\CH}{{\cal H}}
\numberwithin{equation}{section}
\newtheorem{proposition}{Proposition}[section]
\newtheorem{theorem}[proposition]{Theorem}
\newtheorem{corollary}[proposition]{Corollary}
\newtheorem{lemma}[proposition]{Lemma}
\newtheorem{definition}[proposition]{Definition}
\newtheorem{remark}[proposition]{Remark}
\newtheorem{example}[proposition]{Example}
\newtheorem{assumption}[proposition]{Assumption}
\newtheorem{fig}[proposition]{Figure}

\newcommand{\comments}[1]{\marginpar{ \begin{minipage}{0.5in} \tiny #1 \end{minipage} }}
\title{On negative association of some finite point processes on general state spaces}
\author{G\"unter Last
\thanks{Department of Mathematics, Karlsruhe Institute of Technology, Englerstr. 2,  D-76131, Karlsruhe, Germany.}
\\
{\it  Karlsruhe Institute of Technology}
\and Ryszard Szekli
\thanks{Work supported by National Science Centre, Poland grant 2015/19/B/ST1/01152, address: University of Wroc{\l}aw, Mathematical Institute, pl. Grunwaldzki 2/4,  50-384 Wroc{\l}aw, Poland.}
 \\
  {\it  University of Wroc{\l}aw}
} 
\maketitle

\noindent \noindent {\it Abstract:} 
We study negative association (NA) for mixed sampled point processes and show that NA holds for such processes if the random number of points of them fulfills ULC  property. We connect NA property of point processes with dcx  dependence ordering and show some consequences of it for mixed sampled and determinantal point processes. 
Some applications illustrate general theory.

\vspace{0.2cm}

\section{Introduction}

The questions studied in this paper are motivated by several negative dependence properties which are present in combinatorial probability, stochastic processes,  statistical mechanics, reliability and statistics. We focus our study on point processes theory which is a natural tool in many of these fields.  For each of these fields, it seems desirable to get a  better understanding of what it means for a collection of random variables to be {\sl repelling} or mutually negatively dependent.  It is known that it is not  possible to copy the theory of positively dependent random variables.

Early history of various concepts of multivariate negative dependence  are based on topics considered in   Block and Savits \cite{Block1979}, Block, Savits and Shaked \cite{Block1982}, Ebrahimi and Ghosh \cite{ghosh1981multivariate} and Karlin and Rinott \cite{karlin1980classes}. One of the fundamental results discussed in \cite{Block1979} is that if a distribution satisfies an intuitive structural condition called {\sl Condition N} (recalled in the next section), then it satisfies all of the other conditions introduced there.  {\sl Condition N} is satisfied by the multinomial, hypergeometric and Dirichlet distribution as well as several others. It also implies {\sl negative association} introduced by Joag-Dev and Proschan \cite{Joag-Dev1983}. A slightly stronger version of {\sl Condition N} implies a condition based on stochastic ordering (NDS) due to Block, Savits and Shaked \cite{block1985concept}. Negative association has one distinct advantage over the other types of negative dependence. Non-decreasing functions of disjoint sets of negatively associated random variables are also negatively associated. This closure property does not hold for the other types of negative dependence studied in the above mentioned papers. {\sl Condition N} appears in a natural way in the context of queueing networks implying negative association of the population vector in Gordon-Newell networks and negative association of the sojourn times  vector in cyclic networks, see \cite{daduna1996queueing}, and \cite{daduna2004correlation}.

Pemantle in \cite{Pemantle2000} in his negative dependence study confined himself  to binary-valued random variables, in the hope that eliminating the metric and order properties of the real numbers in favour of the two point set, will better reveal what is essential to the questions about negative dependence.
The list of examples that  motivated him to develop techniques for proving that measures have negative dependence properties such as negative association include the uniform random spanning tree, where the vector of indicator functions of the events that the edges of a graph belong to randomly chosen spanning tree is a random vector which is negatively associated, which was proved by   Feder and Mihail \cite{feder1992balanced}. Similar properties hold for weighted spanning trees.  Further items on this list are simple exclusion processes, random cluster models and the occupation status of competing urns. 
Dubhashi and Ranjan \cite{dubhashi1998balls} consider the competing urns example in detail and show negative association of the numbers of balls in each bin and some consequences such as Chernoff bounds for various models (for extensions see \cite{hu2006negative}).
From this fact it follows negative association of the indicators of exceeding any prescribed thresholds in bins. Occupation numbers of urns under various probability schemes have appeared in many places. Instead of multinomial probabilities, one can postulate indistinguishability  of urns or balls and arrive at  Bose-Einstein or other statistics. Negative association  arise in the multinomial models, where Mallows \cite{mallows1968inequality} was one of the first ones to observe negative dependence.

In Borcea et al. \cite{borcea2009negative}  several conjectures related to negative dependence made by Liggett \cite{liggett2002}, Pemantle \cite{Pemantle2000}, and Wagner \cite{wagner2008negatively}, respectively, were solved  and also  Lyons' main results \cite{lyons2003determinantal} on negative association  for determinantal probability measures induced by positive contractions were extended. The authors used several new classes of negatively dependent measures for zero-one valued vectors related to the theory of polynomials and to determinantal measures (for example  strongly Rayleigh measures related to  the notion of proper position for multivariate stable polynomials).
The problem of describing natural negative dependence properties that are preserved by symmetric exclusion evolutions has attracted some attention in the theory of interacting particle systems and Markov processes. In \cite{borcea2009negative} the authors  provide an answer to the aforementioned problem and show that if the initial distribution of a symmetric exclusion process is strongly Rayleigh, then so is the distribution at any time; therefore,  the latter distribution is strongly negatively associated. In particular, this solved an open problem of Pemantle  and Liggett  stating that the distribution of a symmetric exclusion process at any time, with non-random/deterministic initial configuration is negatively associated, and shows that the same is actually true whenever the initial distribution is strongly Rayleigh. In a later paper \cite{liggett2009distributional}, Liggett has applied these results to prove convergence to the normal and Poisson laws for various functionals of the symmetric exclusion process.
Another  possible scenario when utilising negative association is for example : first, show that a model is negatively associated; second, use that negatively associated measures have sub - Gaussian tails; finally use that negative association is known to imply the Chernoff - Hoeffding tail bounds. That kind of approach via strong Rayleigh measures has been realised by Pemantle and Peres \cite{pemantle2014concentration}, and a way of finding negative association via strong Rayleigh property by Peres et al \cite{peres2017random}.

For point processes a negative association result is known in a quite general setting for so called determinantal point processes on locally compact Polish spaces generated by locally trace class positive contractions on natural $L^2$ space, see e.g.\cite{lyons2014determinantal}, Theorem 3.7. A broad list of interesting examples of determinantal point processes can be found in \cite{Soshnikov2000}.
Negative dependence for finite point processes via determinantal and/or strongly Rayleigh measures  have interesting applications in various applied fields such as machine learning, computer vision, computational biology, 
natural language processing, combinatorial bandit learning, neural network compression  and matrix approximations, see for example \cite{anari2016monte}, \cite{kulesza2012determinantal}, \cite{li2015efficient}, \cite{li2016fast}, and references therein.

Another approach to study dependence has been used in finance models. Positive and negative dependence may be seen as some stochastic ordering relation to independence. Such stochastic orderings are called {\sl dependence orderings} (see \cite{joe1997multivariate} or \cite{mullercomparison}). Typical orderings used are supermodular ordering and directionally convex ordering. Relations of these orderings to association and negative association with some applications to concentration inequalities and to the theory of copula functions are given by Christofides and Vaggelatou \cite{Christofides2004}, and Ruechendorf \cite{ruschendorf2004comparison}, respectively. Related results in the theory of point processes and stochastic geometry, where the directionally convex ordering is used to express more clustering in point patterns, are obtained by Blaszczyszyn and Yogeshwaran \cite{blaszczyszyn2015clustering}.

Positive and negative association may be used to obtain information on the distribution of
functionals such as the sum of coordinates. Newman  \cite{newman1984asymptotic} shows that under either a positive or a negative dependence assumption the joint characteristic function of the variables  $\bX$ is well approximated by the product of individual characteristic functions. This allows him to obtain central limit theorems  for stationary
sequences of associated variables. In the positive association case one needs to assume summable covariances, whereas in the negative case one gets this for free. The list of references on central limit theorems for positively/negatively associated variables is very long, a recent reference dealing with point processes is for example the paper by Poinas et al \cite{poinas2017mixing}. 

The central concept of negative dependence in the present paper is negative association ($NA$) which we show for some point processes. The theory and application of $NA$  are not
simply the duals of the theory and application of positive association, but
differ in important respects. Negative association has one distinct advantage
over the other known types of negative dependence. Non-decreasing functions of
disjoint sets of $NA$  random variables are also $NA$. 
Apart from $NA$ property of determinantal point processes not much is known about $NA$ property of other point processes. 
Therefore  it might be interesting to characterise $NA$ property of a different, elementary but very useful class of finite point processes with iid locations of points, so called  mixed sampled point processes.
In order to obtain $NA$ property in this general class of point processes we use some results from the theory of strongly Rayleigh measures on the unit cube (see theorems 3.3 and 3.6). Consequences of $NA$ property of point processes to ordering of dependence in point processes are described in Proposition 4.7 and Corollary 4.8.

\section{Negative association and related definitions}
For the distribution of a real random variable $X$ we say that its density (or probability function) is $PF_2$ if it is log-concave (discrete log-concave) on its support, see eg. Block et al \cite{Block1982} for a detailed description.

The random vector $\bX=(X_1,\ldots, X_n)$ satisfies {\sl Condition N} if there exist $n +1$ independent real random variables $S_0,S_1,\ldots,S_n,$  each having a $PF_2$ density  (or probability function) and a real number $s$ such that
$$\bX\overset{d}{=} [(S_1,\ldots, S_n)|S_0 +S_1 +\cdots +S_n = s]$$
where  $\overset{d}{=}$ denotes equality in distribution, and $[(S_1,\ldots, S_n)|S_0 +S_1 +\cdots +S_n = s]$ denotes a random variable having the distribution of $(S_1,\ldots, S_n)$ conditioned on the event $S_0 +S_1 +\cdots +S_n = s.$

The multinomial distribution is conditional distribution of sums of independent Poisson random variables given that their sum is fixed, and  the Dirichlet is the conditional distribution of independent gammas given that its sum is fixed.

{\sl Condition N} is related to many other definitions of negative dependence, and for example, as summarised in \cite{Block1979}, it is stronger than being completely $RR_2$ in pairs, than having $RR_2$ in pairs measure, than being $S-MRR_2$, and finally stronger than $NA$. We recall the definition of $NA$.

\begin{definition}\label{na}
A random vector ${\bf X}=(X_1,\ldots,X_n)$  is negatively associated  \index{Negatively associated}
(NA)  if, for every subset $A\subseteq \{1,\ldots 
,n\}$
$$
\CV(f(X_{i},i\in A),g(X_{j},j\in A^{c}) ) \le  0,
$$
\noindent whenever $f,g$  are real non-decreasing functions.
\par
\end{definition}
\medskip
\noindent $NA$  may also refer to the set of random variables $\{X_{1},\ldots 
,X_{n}\}$, or to the
underlying distribution of  {\bf X}.

Negative association possesses the following properties (see  Joag-Dev and
Proschan \cite{Joag-Dev1983})
\par
\medskip
\begin{itemize}
\item[(i)] A pair $(X,Y)$  of random variables is $NA$  if and only if
$$
\BP(X\le x,Y\le y) \le  \BP(X\le x) \BP(Y\le y),
$$
\noindent i.e. $(X,Y)$  is negatively quadrant dependent  ($NQD$).
\item[(ii)] For disjoint subsets $A_{1},\ldots 
,A_{m}$  of $\{1,\ldots 
,n\}$, and non-decreasing
positive functions $f_{1},\ldots 
,f_{m},\  {\bf X}$  is $NA$  implies
$$
\EV \prod  ^{m}_{i=1}f_{i}({\bf X}_{A_{i}}) \le  \prod  ^{m}_{i=1} \EV f_{i}({\bf X}_{A_{i}}),
$$
\noindent where $ {\bf X}_{A_{i}} = (X_{j},j\in A_{i}).
$
\item[(iii)] Any (at least two element) subset of $NA$  random variables is $NA$.
\item[(iv)]  If ${\bf X}$  has independent components then it is $NA$.
\item[(v)]  Increasing (non-decreasing) real functions defined on disjoint subsets of a set of $NA$
random variables are $NA$.
\item[(vi)]  If ${\bf X}$  is $NA$ and ${\bf Y}$  is $NA$, and ${\bf X}$ is independent of ${\bf Y}$  then $({\bf X},{\bf Y})$
is $NA$.
\end{itemize}

In some applications negative association appears, when the random
variables are subjected to conditioning.
\begin{theorem}
Let $X_{1},\ldots 
,X_{n}$  be independent, and suppose that
\begin{equation}\label{star1}
\EV (f({\bf X}_{A}) \mid  \sum^{}_{i\in A} X_{i}=s)
\end{equation}
\noindent is increasing in $s$, for every nondecreasing $f$, and every $A\subseteq \{1,\ldots 
,n\}$.
Then the  distribution of $({\bf X \mid } \sum^{n}_{i=1}X_{i}=s)$
is $NA$, for almost all $s$.
\end{theorem}

\par
\medskip
\noindent The above theorem takes on added interest, when considered in conjunction with
the following theorem from \cite{Efron1965}. For a queueing theoretical proof of this theorem see \cite{daduna1996queueing}.
\begin{theorem} 
Let $X_{1},\ldots 
,X_{n}$  be mutually independent with $PF_{2}$  densities and $S_n=\sum^{n}_{i=1}X_{i}.$ Then
$$
\EV (\phi ({\bf X}) \mid S_n=s )
$$
\noindent is increasing in  ( almost every ) $s$, provided $\phi $  is non-decreasing.
\par
\end{theorem}
\begin{corollary}
If $X_{1},\ldots 
,X_{n}$   are independent with $PF_{2}$  densities then the conditional
distribution of $({\bf X \mid } S_n=s )$  is $NA$, for almost all  s.
\par
\end{corollary}
Conditioning with respect to sums is not the only way to obtain $NA$ property by conditioning, as it is shown in \cite{Hu1999}, where conditioning on order statistics were used.

The property of negative association is reasonably useful but hard to verify.
Negatively correlated probability measures appear naturally in many different contexts. Here are some examples related to $NA$ property.
Let ${\bf x}=(x_{1},\ldots 
,x_{n})$  be a set of real numbers. A permutation
distribution is the joint distribution of the vector ${\bf X}$, which takes as
values all permutations of ${\bf x}$  with equal probabilities $1/n$!. Such a
distribution is $NA$. Negatively correlated normal random variables are $NA$. See Joag-Dev and Proschan \cite{Joag-Dev1983}) for these and many other examples.
There are some examples directly related to {\sl Condition N} and conditioning.
$NA$ property of multinomial distributions can be seen from {\sl Condition N}, since it is the conditional distribution of independent
Poisson random variables given their sum. 
One can see that multivariate hyper-geometric distribution is $NA$ because it is the conditional distribution of independent binomial random variables given their sum, see \cite{Block1982}. 
The population vector in Gordon-Newell closed queueing networks
is $NA$, see   \cite{szekli2012stochastic},  Section 3.8, Theorem E.
The sojourn times vector in cyclic queues has got property $NA$, see
\cite{daduna2004correlation}.
A plethora of negative dependence properties and conjectures for zero-one valued vectors  ( $CNA, CNA+, JNRD, JNRD+, h-NLC+, S-MRR_2$, ULC, strongly Rayleigh, Rayleigh, PHR) were introduced and studied in\cite{Pemantle2000} and \cite{borcea2009negativenegative} with an application to symmetric exclusion processes. These classes are related to 
determinantal probability measures described for example in  \cite{lyons2003determinantal} and \cite{lyons2014determinantal}. We shall not recall all of these definitions, and only recall $NA$ property of determinantal point processes.

We shall utilize a slightly broader class than $NA$ in our formulations on dependence orderings. We define this new class of distributions as an analog of weak association in sequence class ($WAS$) introduced by R\"uchendorf in \cite{ruschendorf2004comparison}.  We say that a random vector $\mathbf X$ (or its distribution) is $sNA$ (negatively associated in sequence) if  
\begin{equation}\label{sna}
\CV (\I_{(X_i>t},f(X_{i+1},\ldots ,X_n))\le 0,
\end{equation}
for all $f$ real non-decreasing functions and $t\in \R$, $i=1,\ldots,n-1.$ 

This condition  is equivallent to $[(X_{i+1},\ldots ,X_n)\mid X_i>t)] <_{st} (X_{i+1},\ldots ,X_n)$ for all $t\in \R$ and $i=1,\ldots, n-1$ Here $<_{st}$ denotes the usual strong stochastic ordering on $\R^n.$ For the definitions of strong stochastic orderings we refer to \cite{szekli2012stochastic}.
\subsection{Rayleigh measures and related classes}
In order to study $NA$ property of point processes we  need some other (stronger) properties of negative dependence, which we shall take from \cite{borcea2009negativenegative}. In this section, we will reduce the context of general point processes to point processes which can be represented in distribution by distributions of zero-one valued variables. Therefore we will be first interested  in probability measures on the unit cube. More precisely, let us consider probability function $\mu: 2^{[n]}\to [0,1]$, for $[n]=\{1,\ldots,n\}$, such that $\sum_{S\in 2^{[n]}}\mu (S)=1$. 
For $i\in [n]$  the i-th coordinate function on $2^{[n]}$ is a binary random variable given by $X_i(S) = 1$ if $i \in S$ and $0$ otherwise, where $S \subseteq  [n]$.
There is a 1-1 correspondence between probability functions  on  $2^{[n]}$ and the corresponding generating functions. For $\mu: 2^{[n]}\to [0,1]$ its generating function (multi-affine polynomial), is given by
$$P_{\mu}(\bz) = \sum_{S\subseteq [n]} \mu (S)\bz^S, $$
where $\bz = (z_1,\ldots,z_n),\  S\subseteq[n]$, and $\bz^S=\prod_{i\in S}z_i$.
One can define then that $\mu$ is  $NA$ if 
$$\sum_{S\in 2^{[n]}} f(S)\mu (S)\sum_{S\in 2^{[n]}} g(S)\mu (S)\ge  \sum_{S\in 2^{[n]}} f(S)g(S) \mu(S)$$
for any increasing functions $f, g$ on $ 2^{[n]}$ (with respect to the inclusion) that depend on disjoint sets of coordinates, which is equivalent to the fact that the vector $\bX=(X_1,\ldots,X_n)$  is $NA$ on the probability space $(2^{[n]}, \BP_\mu)$, where $\BP_\mu$ is the probability measure corresponding to probability function (density) $\mu$.
\begin{definition}
 A polynomial with all real coefficients $P \in \R [z_1,\ldots,z_n]$ is called (real) stable if $P(z_1,\ldots, z_ n)= 0$ whenever $ Im (z_j) > 0$,  for $1 \le j \le  n$. 
\end{definition}
The following property implies $NA$.
\begin{definition}
 A measure $\BP_{\mu}$ with density $\mu$ on $2^{[n]}$ (equivalently the corresponding 0-1 vector $\bX$) is called strongly Rayleigh if its generating polynomial $P_\mu$ is (real) stable. 
\end{definition}
\begin{example}\rm
A product measure on $2^{[n]}$ i.e. a measure with density $\mu$ having its  generating function of the form
$$
P_{\mu}(\bz)=\prod_{i=1}^n (p_iz_i+1-p_i),\   0\le p_i \le 1,\ i\in [n],
$$
is strongly Rayleigh.\ $\square$
\end{example}
A useful characterisation of strongly Rayleigh measures comes from the fact that multi-affine polynomial $P\in \R[z_1, \ldots, z_n]$ is (real) stable if and only if
\begin{equation}\label{rayleigh}
\frac{\partial P}{\partial z_i} (\bx) \frac{\partial P}{\partial z_j} (\bx) \ge P(\bx)\frac{\partial ^2}{\partial z_iz_j} P (\bx),
\end{equation}
for all $\bx \in \R^n$, and $1\le i,j\le n$, see \cite{branden2007polynomials}, Theorem 5.6.

A measure $\BP_{\mu}$ with density $\mu$ on $2^{[n]}$ (equivalently the corresponding 0-1 vector $\bX$) is called  Rayleigh if its generating polynomial $P_\mu$ fulfills (\ref{rayleigh}) for all $\bx \in \R^n_+$.

A complex measure $\mu $  on $2^{[n]}$ is called symmetric  if its generating polynomial $P_\mu $ is symmetric in all $n$ variables. The measure $\mu$  is almost symmetric or almost symmetric if $P_\mu$ is symmetric in all but possibly one variable.

One says that $\mu$ is {\it conditionally negatively associated}  ($CNA$) if each measure obtained from $\mu$ by conditioning on some (or none) of the values of the variables is $NA$. Finally, $\mu$ is called {\it strongly conditionally negatively associated} 
$(CNA+)$ if each measure obtained from $\mu$ by imposing external fields and projections is $CNA,$ for details see \cite{borcea2009negative}, where one of the implications in Conjecture 2.6. there, states that  Rayleigh property implies $CNA+$ property of $\mu$.

Pemantle's result (\cite{Pemantle2000}, Theorem 3.7 in §3.5), says that  for symmetric measures, Rayleigh, $CNA$, and $CNA+$ properties are equivalent, and in addition they are equivalent to the fact that the sequence $( \frac{\mu (\sum_{i=0}^nX_i=k)}{\binom{n}{k}}) _{k=0}^n $ is a log-concave sequence. 

A partial answer to the implication from the mentioned Conjecture 2.6 is given in \cite{borcea2009negative} in
Corollary 6.6. where this implication is established for $\mu$  almost symmetric. This implies that (almost) symmetric  Rayleigh measures are $NA$. Strongly Rayleigh measures are $NA$ without assumptions on symmetry. 

Properties for more general polynomials than multi-affine (and distributions of vectors $\bX$ of random variables with values in finite sets of natural numbers) can be reduced to properties of multi-affine polynomials by so called polarization.

Let $P\in  \BC[z_1,\ldots, z_n]$ be a polynomial of degree $d_i$ in the variable $z_i$ for $1\le i \le n$. The polarization $\tilde P$, is the unique polynomial in the variables $z_{ij}, \ 1\le i \le n, 1\le  j \le d_i,$ satisfying
\begin{enumerate}[(1)]
\item $\tilde P$ is multi-affine, 
\item $\tilde P$ is symmetric in the variables $z_{i1},\ldots ,z_{id_i},$
for $1\le i \le n,$
\item
 if we let $z_{ij} =z_i $  for all $ i, j$  in $\tilde P ,$  we recover $P$ .
\end{enumerate}
It is known (see e.g. \cite{borcea2009negative}, Corrolary 4.7) that $P\in \BC[z_1,\ldots, z_n]$ is stable iff $\tilde P$ is stable.
More precisely, if $( {X_{i,j} : i \le n, j\le  d_i})$ is a finite family of non-negative integer variables with real stable generating polynomial, the aggregate variables $X_i=\sum_{j=1}^{d_i} X_{i,j}$
 will also have a real stable generating polynomial. This is because it follows from the definition that stability is preserved by substituting $z_{i,j} = z_i$ for all $i,j$. Conversely, if $(X_i : 1 \le  i \le n)$ are random variables whose joint law has a real stable generating polynomial $P$, one can define the polarization of $P$ by the substitutions $ z_i^j=\binom{d_i}{j}^{-1} e_j(z_{i,1}, \dots, z_{i,d_i} )$ where $d_i$ are upper bounds for the values of $X_i$  and $e_j$  is the elementary
symmetric function of degree $j$  on $d_i$ variables. Stability of $P$ implies stability of the polarization $\tilde P$  of $P$, hence the strong Rayleigh property of a collection of binary variables $(X_{i,j})$ with probability generating function $\tilde P$.

\section{NA for mixed sampled point processes}

We shall adopt our  notation from the book by Last and Penrose \cite{last2017lectures}. Let $(\Omega, \cal F, \mathbb P)$ be a probability space and  $(\mathbb X , {\mathcal X})$ a measurable state space.
Denote by $ \mathbf N_{<\infty}(\mathbb X)=\mathbf N_{<\infty}$  the space of all measures $\mu$ on $(\mathbb X , {\mathcal X})$ such that $\mu(B) \in \Z_+$ for all $B \in {\mathcal X}.$ Let $\mathbf N(\mathbb X)=\mathbf N$ be the space of all measures that can be written as a countable sum of measures
from $\mathbf N_{<\infty}$. An example is the Dirac measure $\delta _x$ for a point $x \in \mathbb X ,$ given by $\delta_x(B):=\I_B (x)$. 

We define a {\sl point process} $\eta $ as a measurable mapping from $(\Omega, \cal F, \mathbb P)$ to $(\mathbf N, {\mathcal N})$ (${\mathcal N}$ is the smallest $\sigma$-field on $\mathbf N$ such that $\mu \mapsto \mu(B)$ is measurable for all $B \in {\mathcal X}) $.

We restrict our attention in this paper to  finite point process with points located in a complete separable metric space $\BX$ with ${\mathcal X}$ being the Borel $\sigma$ field.

We define NA property of point processes as follows.
\begin{definition}\label{na-pp}
A point process $\eta$ is negatively associated ({NA}) or negatively associated in sequence ($sNA$)  if for each collection of disjoint sets $B_1,\ldots, B_n\in {\mathcal X}$ the vector $(\eta (B_1),\ldots,\eta (B_n))$ is $NA$ or $sNA$, respectively, as defined for random vectors.
\end{definition}
For a more general definition of negative association for point processes (random measures) we refer to \cite{lyons2014determinantal} or \cite{yogesh-na}. 

For a Borel set $A \subseteq \mathbb X$, let ${\mathcal N}_{A}$ denote the $\sigma$-field on $\mathbf N$ generated by the functions
$\mu \mapsto \mu (B)$  for Borel $B \subseteq A$.  The natural (inclusion) partial order on $\mathbf N$ allows us to define $ f: \mathbf N\to \R$ which is increasing.  We say that a point process $\eta$  has {\it negative associations} if $\BE (f (\eta)g(\eta))\le \BE (f(\eta))\BE(g(\eta))$ for every pair $f, g$ of real bounded increasing functions that are measurable with respect to complementary subsets $A$, $A^c$
of $\mathbb X$, meaning that a function  is measurable with respect to $A$ if it is measurable with respect to ${\mathcal N}_{A}$. It is clear that with such a general definition we have that if $\eta$ has {\it negative associations} then $\eta$ is $NA$. Let us recall Theorem 3.7 from \cite{lyons2014determinantal}. Let $\lambda$ be a Radon measure on a locally compact Polish space $\mathbb X$. Let $K$ be a locally trace-class positive contraction on $L_2(\mathbb X, \lambda)$. By $\eta_K$ we denote the determinantal point process generated by $K$, for details see \cite{lyons2014determinantal}, section 3.2.
\begin{theorem}\label{determinantal-na}
 The determinantal point process $\eta_K$ defined above has negative associations.
\end{theorem}

Apart from determinantal point processes not much is known about $NA$ property of point processes. 
Therefore  we concentrate our efforts  on characterising $NA$ property for an elementary but very useful class of finite point processes with iid locations of points, more precisely
our main focus in this paper is on the  class of so  called {\sl  mixed sampled point processes} on $\mathbb X$,  defined by 
\begin{equation}\label{mixed}
\eta =\sum_{i=1}^\tau \delta_{X_i},
\end{equation}
where $(X)_{i\ge 1}$ is iid with distribution $F$, and $\tau\in \N\cup \{0\}$ is independent of $(X_i)_{i\ge 1}$.

For this process, given any finite partition $A_1,\ldots,A_k$ of $\mathbb X$, conditionally on $\tau$, the joint distribution of the number of points is given by
$$
\BP(\eta (A_1)=n_1,\ldots,\eta (A_k)=n_k|\tau =N)={{N}\choose {n_1\cdots n_k} }F(A_1)^{n_1}\cdots F(A_k)^{n_k},
$$
and unconditionally
$$
\BP (\eta (A_1)=n_1,\ldots,\eta (A_k)=n_k)=\sum_{N=0}^{\infty} \BP (\tau=N){{N}\choose {n_1\cdots n_k} }F(A_1)^{n_1}\cdots F(A_k)^{n_k}.
$$
The joint probability generating function is therefore given by
\begin{equation}\label{pgf}
\EV (z_1^{\eta(A_1)}\cdots z_k^{\eta(A_k)})=P _{\tau}(F(A_1)z_1+\cdots+F(A_k)z_k),
\end{equation}
where $P _\tau(z)=\EV (z^\tau)$, $z\in [0,1]$.

First, we shall consider mixed point processes defined by (\ref{mixed}) for which random variables $\tau$ are of the form $$\tau=\sum_{i=1}^nU_i$$
where $n\in\N$ and 
$U_1,\ldots, U_i $ are independent Bernoulli variables with possibly different success probabilities. The class of random variables which are the sums of $n$ Bernoulli variables we denote by 
$
{\mathcal Q}_n.
$
Moreover we denote by
$$
{\mathcal Q}:= \cl (\cup_{n=1}^\infty\BQ_n)
$$
the class of all distributions with supports contained in $\{0,1,\ldots\}$ appearing as  weak limits of distributions from $\BQ_n$, $n\ge 1,$ i.e. the weak closure of $\cup_{n=1}^\infty\BQ_n$ 
The main results of this paper are  contained in Theorems \ref{finiteNA1} and \ref{finiteNA2}.
\begin{theorem}\label{finiteNA1}
Suppose that $\eta$ is a mixed sampled point process on $\mathbb X$, defined by (\ref{mixed}), for which $\tau\in {\mathcal Q}$. Then $\eta$ is $NA$.
\end{theorem}
\proof 
Let $B_1,\ldots,B_n \in {\mathcal X}$ be a partition of $\BX$, and $q_i:=F(B_i)$, $i=1,\ldots,m$. Define by
$$
\Z _i:=(\I_{\{X_i\in B_1\}},\ldots , \I_{\{X_i\in B_m\}},
$$
the vector generated by the $i$-th sample $X_i\in \BX$, $i\ge 1$. Note that each $\Z_i$ has multinomial distribution with success parameters $q_1,\ldots,q_m$ and the number of trials equal 1, and as such is $NA$. Moreover $\Z_i, i\ge 1$ are independent.  
Let $\bU=(U_1,\ldots, U_n)$ be a vector of zero-one valued, independent random variables which is independent of $\Z_i, i\ge 1.$ The vector composed as $(\bU, \Z_1,\ldots, \Z_n)$ is $NA$ because of properties (iv) and (vi) of $NA$.

 Now using property (v) we get that the vector $(U_1\Z_1,\ldots,U_n\Z_n)$ is $NA$ as a monotone transformation (multiplication)  of disjoint coordinates of $(\bU, \Z_1,\ldots, \Z_n).$ Again using property (v), this time to $(U_1\Z_1,\ldots,U_n\Z_n)$ and using appropriately addition we get that the vector $\sum_{i=1}^n U_i\Z_i$ is $NA.$ It is clear that $\sum_{i=1}^n U_i\Z_i$ has got the same distribution as $\sum_{i=1}^{\tau} \Z_i$ if $\tau \in {\mathcal Q}_n$, which in turn is the same as for $(\eta (B_1),\ldots, \eta (B_m)).$ This finishes the proof for $\tau \in {\mathcal Q}_n$ for arbitrary $n\in \N.$ 

For $\tau \in \mathcal Q$ there exist a sequence $\tau_k\overset{d}{\longrightarrow}\tau,$  $k\to \infty,$ for $\tau_k\in \cup_{n=1}^\infty\mathcal Q_n$, and
$$
\BE \bigg[ f\big(\sum_{i=1}^{\tau_k}\Z_i\big) g\big(\sum_{i=1}^{\tau_k}\Z_i\big)\bigg]\le \BE \bigg[f\big(\sum_{i=1}^{\tau_k}\Z_i\big)\bigg]\BE\bigg[g\big(\sum_{i=1}^{\tau_k}\Z_i\big)\bigg],
$$
for $f,g$ supported by disjoint coordinates, which are non-decreasing and bounded. Letting $k\to \infty$ gives 
$$
\BE \bigg[f\big(\sum_{i=1}^{\tau}\Z_i\big)g\big(\sum_{i=1}^{\tau}\Z_i\big)\bigg]\le \BE \bigg[f\big(\sum_{i=1}^{\tau}\Z_i\big)\bigg]\BE\bigg[g\big(\sum_{i=1}^{\tau}\Z_i\big)\bigg].
$$
Since each non-decreasing function can be monotonically approximated by non-decreasing and bounded functions, we get $NA$ property of $\eta$. $\square$

The class $\mathcal Q$ can be completely characterized, see e.g. \cite{aleman2004real}.
\begin{lemma}
$$\tau\in \mathcal Q\ \   {\it iff} \ \  \tau=^d \tau_1 +\tau_2,$$
where $\tau_1, \tau_2$ are independent and $\tau_1$ has Poisson distribution and $\tau_2=^d \sum _{i=1}^\infty U_i$ for independent zero-one valued variables $U_i$ with $p_i=\BP(U_i=1)\ge 0$, $i\ge 1,$ such that $\sum_{i=1}^\infty p_i<\infty.$
\end{lemma}
It is interesting to note that hypergeometric random variables belong to the class $\mathcal Q$, see e.g. \cite{hui2014representation}

We say that a  real sequence $(a_i)_{i=0}^n$ has no internal zeros if the indices of its non-zero terms form a discrete interval. Following Pemantle \cite{Pemantle2000} we shall use the following class of sequences and distributions.
\begin{definition}
W say that a finite real sequence $(a_i)_{i=0}^n$ of non-negative real numbers with no internal zeros is {\it ultra log-concave} {\rm (ULC(n))} if $$
\left(\frac{a_i}{{{n}\choose {i}}}\right)^2\ge \frac{a_{i-1}}{{{n}\choose {i-1}}}\frac{a_{i+1}}{{{n}\choose { i+1}}}, \ \ i=1,\ldots,n-1.
$$
\end{definition}
The class of random variables for which their probability functions have the above property we denote by $\mathcal S_n$, i.e.
$$
\mathcal S_n:=\{\tau : (\BP(\tau=i))_{i=0}^n\ \  {\rm is} \ \ {\rm ULC(n)}\}.
$$
It is known that if a non-nengative sequence $(a_i)_{i=0}^n$ is ULC(n), and a nonnegative sequence $(b_i)_{i=0}^m$ is ULC(m) then the convolution of these sequences is ULC(m+n), see \cite{liggett1997ultra}, Theorem 2. Let
$$
\mathcal S:=\cl (\cup_{n=1}^\infty\mathcal S_n).
$$
Sums of independent variables from the class $\mathcal S$ are in $\mathcal S$. We shall see below that $\mathcal Q\subseteq \mathcal S$. Utilizing the class $\mathcal S$, the Theorem \ref{finiteNA1} can be generalised with a use of elementary symmetric functions.
\begin{theorem}\label{finiteNA2}
Suppose that $\eta$ is a mixed sampled point process on $\mathbb X$, defined by (\ref{mixed}), for which $\tau\in \mathcal S$. Then $\eta$  is $NA$.
\end{theorem}
\proof 
Assume first that $\tau\in \mathcal S_n$.
Let $B_1,\ldots,B_n \in {\mathcal X}$ be a partition of $\BX$, and $q_i:=F(B_i)$, $i=1,\ldots,m$. Define by
$$
\Z _i:=(\I_{\{X_i\in B_1\}},\ldots , \I_{\{X_i\in B_m\}},
$$
the zero-one valued vector generated by the $i$-th sample $X_i\in \BX$, $i\ge 1$. Note that each $\Z_i$ has multinomial distribution with success parameters $q_1,\ldots,q_m,$ and the number of trials equal 1, and as such is $NA$. Moreover $\Z_1,\ldots$ are independent.  

For fixed $n\in \N$,  let $\bU=(U_1,\ldots, U_n)$ be, independent of $\Z_i, i\ge 1,$  vector,  of 0-1 valued random variables obtained in the following way.
For the generating function of $\tau$, $P_{\tau}(z)=\BE(z^\tau)$, we define $\bU$ by providing its multidimensional generating function. It is obtained by substituting in $P_{\tau}(z)$, for each $k=0,\ldots, n$, 
$$z^k:=\binom{n}{k}^{-1} e_k(z_{1}, \dots, z_{n} ),$$
 where $e_k(z_{1}, \dots, z_{n} )$ is  the k-th elementary symmetric polynomial. Note that immediately from the definition of the elementary symmetric polynomials, for each $k$, the function $\binom{n}{k}^{-1} e_k(z_{1}, \dots, z_{n} )$ of variables $z_{1}, \dots, z_{n} $ is the multivariate generating function of a vector of $n$ 0-1 valued variables which takes on exactly $k$ values 1 with the same probability $\binom{n}{k}^{-1}$ defined for all possible selections of $k$  coordinates on which the values 1 are obtained. The distribution of $\bU$ defined in such a way is the mixture with the coefficients $a_k:=\BP(\tau=k)$ of the distributions corresponding to $\binom{n}{k}^{-1} e_k(z_{1}, \dots, z_{n} )$, $k=0,\ldots, n$. Since each function $e_k$ is symmetric in variables $z_{1}, \dots, z_{n} $, the same is true for the generating function of $\bU=(U_1,\ldots, U_n)$, therefore  $(U_1,\ldots, U_n)$ are exchangeable. Moreover  $\sum_{i=1}^nU_i\overset{d}{=} \tau$, since by setting $z_1=\cdots=z_n:=z$ we obtain $P_{\tau}(z)$. In other words  the sequence $(\BP(\tau=i))_{i=0}^n$ is the rank sequence for the vector $\bU=(U_1,\ldots, U_n)$.

From our assumption  the rank sequence for $\bU=(U_1,\ldots, U_n)$ is ULC(n) and from Theorem 2.7 in  \cite{Pemantle2000} we obtain that $\bU$ is $NA$. Now the vector composed as $(\bU, \Z_1,\ldots, \Z_n)$ is $NA$ because of property (vi) of $NA$. Using property (v) we get that the vector $(U_1\Z_1,\ldots,U_n\Z_n)$ is $NA$ as a monotone transformation (multiplication)  of disjoint coordinates of $(\bU, \Z_1,\ldots, \Z_n).$ Again using property (v), this time to $(U_1\Z_1,\ldots,U_n\Z_n)$ and using appropriately addition we get that the vector $\sum_{i=1}^n U_i\Z_i$ is $NA.$ It is clear that $\sum_{i=1}^n U_i\Z_i$ has got the same distribution as $\sum_{i=1}^{\tau} \Z_i$, which in turn has the same  distribution as $(\eta (B_1),\ldots, \eta (B_m)).$ This finishes the proof for $\tau \in \mathcal S_n$, for arbitrary $n\in \N.$ For $\tau \in \mathcal S$, we apply an analogous limiting argument as in Theorem \ref{finiteNA1}. $\square$

The following lemma may be regarded as known since it is an immediate consequence of the classical Newton inequalities.  We put its formulation in the setting of the introduced in this paper classes of random variables. 
\begin{lemma}\label{newton}
For all $n>1$
$$
\mathcal Q_n \subseteq \mathcal S_n.
$$
\end{lemma}
\proof
Suppose $\tau\in\mathcal Q_n$. Then, for its generating function, 
$$P_\tau(z)=(1-p_1+p_1z)\cdots (1-p_n+p_nz)=p_1\cdots p_n(\frac{1-p_1}{p_1}+z)\cdots (\frac{1-p_n}{p_n}+z).$$
For $a_k:=\frac{1-p_k}{p_k}$ we have
$$
P_\tau(z)=p_1\cdots p_n(a_1+z)\cdots (a_n+z)=p_1\cdots p_n[x^n+c_1x^{n-1}+\cdots+c_n],
$$
where $c_n=a_1+\cdots+a_n$, $c_2=a_1a_2+\cdots+ a_{n-1}a_n$, $\cdots$, $c_n=a_1\cdots a_n,$ i.e. $c_k$, $k=0,\ldots, n$, coefficients are given by the corresponding elementary symmetric polynomials in variables $a_i, i=1,\ldots, n$. It is known from the classical Newton's inequalities that for $k=1,\ldots,n-1$
$$
\left( \frac{c_k}{\binom{n}{k}}\right) ^2\ge  \frac{c_{k-1}}{\binom{n}{k-1}}\frac{c_{k+1}}{\binom{n}{k+1}},
$$
and since $\BP(\tau=k)=p_1\cdots p_n c_{n-k}$, we get that the sequence $(\BP(\tau=k))_{k=0}^n\ \  {\rm is} \ \ $ULC(n), and therefore $\tau \in\mathcal S_n$. $\square$ 

Lemma \ref{newton} implies the announced before inclusion  $\mathcal Q\subset \mathcal S$.

It is interesting to note that arguments utilised in  Theorem \ref{finiteNA2} can be used for random vectors with arbitrary positive values.

\begin{proposition}
Assume that $\Z_i=(Z_i^1,\ldots,Z_i^m)$, $i\ge 1$ is a sequence of independent, identically distributed  random vectors
with components in $\R_+$ such  that for each $i\ge 1$
$\sum^m_{j=1} \I\{Z_i^j>0\}\le 1$, that is, at most one of the
components can be positive. Then for $\tau \in \mathcal S$ which is independent of $\Z_i, \ i\ge 1$, the vector $ \bW:=\sum^\tau_{i=1}\Z_i$ is $NA$.
\end{proposition}
\proof We use basically the same argument as in Theorem \ref{finiteNA2}. Let $\bU=(U_1,\ldots, U_n)$ be, independent of $\Z_i, i\ge 1$, vector of 0-1 valued random variables obtained by its generating function obtained as follows. In  $P_{\tau}(z)$, we substitute $z^k:=\binom{n}{k}^{-1} e_k(z_{1}, \dots, z_{n} )$, $k=1,\ldots,n$, where $e_k(z_{1}, \dots, z_{n} )$ are the elementary symmetric polynomials. This substitution defines a generating function of variables $z_1,\ldots,z_n$. It is then immediate  that $\sum_{i=1}^nU_i=^d \tau$, that is the sequence $(\BP(\tau=i))_{i=0}^n$ is the rank sequence for (symmetric) $\bU=(U_1,\ldots, U_n)$. From our assumption  we have that the rank sequence for $\bU=(U_1,\ldots, U_n)$ is ULC(n) and from Theorem 2.7 in  \cite{Pemantle2000} we obtain that $\bU$ is $NA$. Now the vector composed as $(\bU, \Z_1,\ldots, \Z_n)$ is $NA$ because of the following lemma and property (vi) of $NA$.
\begin{lemma}\label{na1}
Assume that $ \Z=(Z^1,\ldots,Z^m)$ is a random vector
with components in $\R_+$. Assume that
$\sum^m_{j=1} \I\{Z^i>0\}\le 1$, that is, at most one of the
components can be positive. Then $\Z$ is $NA$.
\end{lemma}
{\it Proof of lemma} \ref{na1}. In order to show that
$\CV (f(Z^1,\ldots,Z^k),g(Z^{k+1},\ldots,Z^m))\ge 0$
for non-decreasing $f$ and $g$, it suffices to assume that $f(0)=g(0)=0$. Otherwise one can consider $f-f(0)$ and $g-g(0)$. Because only of the coordinates can be non-zero, we get 
$E[f(Z^1,\ldots,Z^k)g(Z^{k+1},\ldots,Z^m))]=0$,
while the  product of the expectations is non-negative
since, $ f\ge 0$ and $g\ge 0$. $\square$

Now, using property (v),  we get that the vector $(U_1\Z_1,\ldots,U_n\Z_n)$ is $NA$ as a monotone transformation (multiplication)  of disjoint coordinates of $(\bU, \Z_1,\ldots, \Z_n).$ Again using property (v), this time to $(U_1\Z_1,\ldots,U_n\Z_n)$ and using appropriately addition we get that the vector $\sum_{i=1}^n U_i\Z_i$ is $NA.$ It is clear that $\sum_{i=1}^n U_i\Z_i$ has got the same distribution as $\sum_{i=1}^{\tau} \Z_i$. $\square$

The above proposition can be used to study random measures other than point processes. We shall pursue this topic elsewhere. 

\section{Dependence orderings  for point processes}

An extensive study of dependence orderings for multivariate point processes on $\R$ is contained in \cite{kulik2005dependence}. Related results in the theory of point processes and stochastic geometry, where the directionally convex ordering is used to express more clustering in point patterns, are obtained by Blaszczyszyn and Yogeshwaran \cite{blaszczyszyn2015clustering}, see also references therein.
We shall use $NA$ property of point processes we study to obtain comparisons related to dependency properties. First we recall some  basic facts on dependence orderings of vectors and their relation to $NA$ which can be directly utilized for point processes.

\subsection{Dependence orderings and negative correlations for vectors}

For a function $f:\R^n\to\R $ define the difference operator $\Delta^\epsilon_i,$ $\epsilon >0,$ $1\le i\le n$
by
$$\Delta^\epsilon_i f(\bx) = f(\bx+\epsilon \be_i)-f(\bx)$$
where $\be_i$ is the i-th unit vector. Then $f$ is called
{\it super-modular} if for all $1\le  i < j \le n$ and $\epsilon, \delta > 0$ $$\Delta^\delta_j\Delta^\epsilon_i f(\bx)\ge 0.$$ for all $\bx\in \R ,$ and {\it directionally convex} if this inequality holds for all $1\le i \le j \le n.$
Let $\CF^{sm}, \CF^{dcx}$ denote the classes of super-modular, and directionally convex functions.
 Then of course $\CF^{dcx}\subseteq \CF^{sm}.$ 
Typical examples from  $\CF^{dcx}$ class of functions are $ f(\bx)=\psi (\sum_{i=1}^nx_i)$, for $\psi$ convex, or $f(\bx)=\max_{1\le i\le n}x_i$, but there are many other useful functions in this class, see for example \cite{blaszczyszyn2009directionally}.

The corresponding stochastic orderings one defines by
$\bX <_{sm}\bY$ if $\BE f(\bX)\le \BE f(\bY)$ for all $f\in \CF^{sm},$ and analogously for $\bX<_{dcx}\bY.$ For differentiable functions $f$ one obtains $ f \in \CF^{sm}$ iff $\frac{\partial ^2}{\partial x_ix_j}f\ge 0,$ for $i < j,$  and $ f \in \CF^{dcx}$ iff this inequality holds for $i \le j.$ While comparison of $\bX$ and $\bY$  with respect to $<_{sm}$  implies (and is restricted to the case of)  identical marginals $X_i=^dY_i,$ the comparison with respect to the smaller class $\CF^{dcx}$ implies  convexly increasing marginals $X_i <_{cx} Y_i$ (which means by definition that $\BE \psi(X_i)\le \BE\psi (Y_i)$ for all $\psi:\R\to\R$ convex). Both of these orderings belong to the class of so called {\it dependency orderings}, see e.g. \cite{joe1997multivariate}, which is defined by a list of suitable properties, among them the property that $\CV(X_i,X_j)\le \CV (Y_i,Y_j)$. 

In \cite{ruschendorf2004comparison} another stochastic ordering related to dependence comparisons was introduced by the condition $\bX<_{wcs}\bY$ ({\it weakly conditional increasing in sequence} order) iff 
$$
\CV (\I_{(X_i>t},f(X_{i+1},\ldots ,X_n))\le \CV (\I_{(Y_i>t}, f(Y_{i+1},\ldots ,Y_n)),
$$
for all $ f $ monotonically non-decreasing, and all $t\in \R$, $1\le i\le n-1.$
The following theorem from \cite{ruschendorf2004comparison} connects the above defined orderings.
\begin{theorem}\label{wcs}
Let $\bX, \bY $ be $n$-dimensional random vectors.
\begin{enumerate}
\item If $X_i\overset{d}{=}Y_i, \ 1\le i\le n$ then $\bX<_{wcs}\bY$ implies that $\bX<_{sm}\bY.$
\item If $X_i<_{cx}Y_i,\ 1\le i\le n$ then $\bX<_{wcs}\bY$ implies that $\bX<_{dcx}\bY.$
\end{enumerate} 
\end{theorem}
Dependence orderings can be used to define some classes of distributions with negative or positive covariances when applied to vectors with independent components. More precisely, denote by $\bX^*$  a vector with independent components, and $X_i^*\overset{d}{=} X_i$.
Then
$\bX $ is called {\it weakly associated in sequence}  if $\bX^* <_{wcs} \bX$.
Note that $\bX$ is $sNA$, as defined by (\ref{sna}), iff $\bX<_{wcs}\bX^*.$
Note also  that  $sNA$ is further equivalent to
$$[(X_{i+1},\ldots ,X_n)\mid X_i>t)] <_{st} (X_{i+1},\ldots ,X_n) ,$$ 
for all $i=1,\ldots, n-1$, $t>0$, where $[(X_{i+1},\ldots ,X_n)\mid X_i>t)]$ denotes a random vector which has got the distribution of $(X_{i+1},\ldots ,X_n)$ conditioned on the event $\{X_i>t\}$, and $<_{st}$ is the usual (strong) stochastic order.
It is clear directly from the definition that $NA$ property implies $sNA$ property. For $\bX$ being $sNA$, Theorem  \ref{wcs} implies  that, for example, (see also \cite{Christofides2004} for the case $NA$)
$
\sum_{i=1}^n X_i <_{cx}\sum_{i=1}^n X^*_i,
$
and
$
\max_{1\le k\le n}\sum_{i=1}^k X_i<_{icx}\max_{1\le k\le n}\sum_{i=1}^k X^*_i,
$
 where $<_{icx}$ is defined similarly as $<_{cx}$ but with the use of non-decreasing convex functions.
Taking other super-modular functions it is possible to get maximal inequalities for $sNA$ vectors as in \cite{Christofides2004} . For $<_{dcx}$ ordering we get the following corollary from Theorem \ref{wcs} which will be used for point processes with $NA$ property.
\begin{corollary}\label{na-dcx}
Suppose $\bX$ is $sNA$ and $\bY^*$ has independent coordinates with $X_i<_{cx}Y_i^*$ then $\bX<_{dcx}\bY^*$.
\end{corollary}

\subsection{NA and dependence orderings for point processes}

We shall connect the approach by dependence orderings with negative association and by using Corollary \ref{na-dcx} we shall be able to compare dependence for  point processes which posses $NA$ property with some other processes. This will   imply in particular comparison results for void probabilities and/or moment measures.


\begin{definition}
Two point processes $\eta_1,\eta_2$ on a complete, separable metric space $\mathbb X$ are ordered in directionaly convex order (weakly conditional increasing sequence order) $$\eta_1<_{dcx}\eta_2\   (\eta_1<_{wcs}\eta_2)\ \ {\rm iff}  \ (\eta_1(B_1),\ldots, \eta_1 \ (B_n))<_{dcx}(<_{wcs})\  (\eta_2(B_1),\ldots, \eta_2 (B_n)),$$ 
as defined for random vectors, for all disjoint, bounded  Borel sets $B_1,\ldots, B_n$, $n\ge 1.$
\end{definition}

The following consequence of $<_{dcx}$ ordering for point processes is known from \cite{blaszczyszyn2015clustering}, Proposition 6.
\begin{lemma}
Let $\eta_1,\eta_2$ be two point process on $\R^d$. If $\eta_1<_{dcx}\eta_2$ then
\begin{enumerate}
\item (moment measures) $\BE (\eta_1(B_1)\cdots\eta_1 (B_n)))\le \BE (\eta_2(B_1)\cdots\eta_2 (B_n)))$ \\ for all disjoint, bounded Borel sets $B_1,\ldots, B_n$,\item (void probabilities) $\BP(\eta_1(B)=0)\le \BP(\eta_2 (B)=0)$ \\ for all bounded Borel sets $B$.
\end{enumerate}
\end{lemma}
Using $<_{wcs}$ criterion for $<_{dcx}$ from Theorem \ref{wcs}, we get
\begin{corollary}
Let $\eta_1,\eta_2$ be two point process on $\R^d$. If $\eta_1(B)<_{cx}\eta_2(B)$ and $\eta_1<_{wcs}\eta_2$ then the moment measures and void probabilities comparisons from the above lemma hold.
\end{corollary}
An interesting case for such comparisons is when $\eta_2$ is a Poisson point process. 
\begin{proposition}\label{momentrd}
Suppose $\eta$ is a simple point process on $\R^d$ which is {\rm sNA}  then 
\begin{enumerate}
\item (moment measures) $\BE (\eta (B_1)\cdots\eta  (B_n)))\le \BE (\eta (B_1))\cdots\BE(\eta (B_n))$ \\ for all disjoint, bounded Borel sets $B_1,\ldots, B_n$, which implies $$ \BE (\exp(\int_{\R^d}h(\bx)\eta(d\bx)\le \exp(\int_{\R^d}(e^{h(\bx)}-1)\BE\eta(d\bx)),$$
for all $h\ge 0$.
\item (void probabilities) $\BP(\eta (B)=0)\le \exp(-\BE\eta(B))$\\  for all bounded Borel sets $B$, which implies 
\begin{equation*} \label{void}\BE (\exp(-\int_{\R^d}h(\bx)\eta(d\bx)\le \exp(\int_{\R^d}(e^{-h(\bx)}-1)\BE\eta(d\bx)),\end{equation*}
for all $h\ge 0$.
\end{enumerate}
\end{proposition}
\proof From our assumption $(\eta (B_1)\cdots\eta  (B_n))$ is $sNA$. Since $\bX$ is $sNA$, as defined by (\ref{sna}) iff $\bX<_{wcs}\bX^*$ then  from  Theorem \ref{wcs}, $(\eta (B_1)\cdots\eta  (B_n))<_{sm}\bX^*$, where $\bX^*$ denotes the corresponding independent version. From the definition of $<_{sm}$ we obtain directly that $\BE (\eta (B_1)\cdots\eta  (B_n)))\le \BE (\eta (B_1))\cdots\BE(\eta (B_n)).$ Now from Proposition 1 in \cite{blaszczyszyn2015clustering} we obtain that this implies $ \BE (\exp(\int_{\R^d}h(\bx)\eta(d\bx)\le \exp(\int_{\R^d}(e^{h(\bx)}-1)\BE\eta(d\bx)),$ for all $h\ge 0$.
 
\noindent Regarding void probabilities, since $(\eta (B_1)\cdots\eta  (B_n))$ is $sNA$ for all  bounded Borel sets $B_1,\ldots, B_n,$ it follows that $(\eta(B),\eta(B'))$ is also $sNA$ for all bounded Borel $B,B'$. Directly from the definition of $sNA$ we conclude that $\BP(\eta (B)=0),\eta (B')=0)\le \BP(\eta (B)=0)\BP(\eta (B')=0)$. Now, from Proposition 3.1 in \cite{blaszczyszyn2014comparison} we get that $\BP(\eta (B)=0)\le \exp(-\BE\eta(B)),$ for all bounded Borel sets $B$. Moreover, from Proposition 2 in \cite{blaszczyszyn2015clustering} this inequality  holds if and only if $\BE (\exp(-\int_{\R^d}h(\bx)\eta(d\bx)\le \exp(\int_{\R^d}(e^{-h(\bx)}-1)\BE\eta(d\bx)),$
for all $h\ge 0$.
$\square$

From Corollary \ref{na-dcx}, the above result can be  modified (using the same argument)  to obtain  $<_{dcx}$ ordering for processes on general state spaces.  For point processes on $\R^d$, this result can also be used to obtain  comparison of moment measures and void probabilities for $sNA$ processes with Poisson processes.
\begin{proposition}\label{dcx-poi}
Suppose $\eta_1$ is a  point process on  a complete, separable metric space $\mathbb X$  which is {sNA},  and $\eta_1 (B)<_{cx} Y$, for all bounded Borel sets $B$, where $Y$ is a random variable with Poisson distribution $Po(\BE\eta_1(B))$. Then $$\eta_1<_{dcx}\eta_2,$$ where $\eta_2$ is a Poisson point process with intensity measure $\BE\eta_1.$ 
\end{proposition}

\subsection{$<_{dcx}$ comparisons for mixed sampled and determinantal point processes}

Before we formulate  more general results, which use our results on $NA$ properties of mixed sampled point processes,  we start with an example which gives a direct approach to a comparison of void probabilities.

\begin{example}[Comparison of binomial mixed sampled p.p. with Poisson p.p. on $\R^d$]\rm
Let $\eta$ be mixed sampled point process on $\R^d$ with $\tau$ being binomially distributed random variable. We shall compare void probabilities for this process with void probabilities of Poisson point process with the same intensity measure.  In general,
for a simple point processes $\eta$ on $\R^d$ in order to get $\BP(\eta (B)=0)\le \exp(-\BE\eta(B))$, it is enough to check (see Proposition 3.1, \cite{blaszczyszyn2015clustering}) whether
$$
\BP(\eta(B)=0, \eta(B')=0)\le \BP(\eta(B)=0))\BP(\eta(B')=0),
$$
for disjoint $B, B'.$
We have
for arbitrary, measurable disjoint sets $B$, $B'$
\begin{align*}
\BP(\eta (B)=\eta (B')=0)
&=\BP(\tau =0)+\sum_{n=1}^\infty \BP(\tau =n) \BP(X_1\notin B\cup B', \ldots . X_n\notin B\cup B')\\
&= \BP(\tau =0)+\sum_{n=1}^\infty \BP(\tau =n)(1-F(B\cup B'))^n\\
&=P_\tau (1-F(B\cup B')) =P_\tau (1-(F(B)+F( B')) .
\end{align*}

 Since $\tau$ has a binomial distribution with, say,  parameters $n$ (number of trials), and $p\in (0,1)$ (success probability) then
$P_\tau (1-s)=(p(1-s)+(1-p))^n$. It is easy to see by differentiaition that $\phi (s):=-\log P_\tau (1-s)$ is then an increasing and convex function such that $\phi(0)=0$.
It is known that such a function is superadditive (see e.g. Bruckner and Ostrow (1963)), therefore $\phi(s+t)\ge \phi(s)+\phi (t)$, and then
$P_\tau(1-(F(B)+F(B'))\le P_\tau (1-F(B))P_\tau (1-F(B'))$. In this case, for disjoint $B, B'$ we obtain
$$\BP(\eta (B)=\eta (B')=0)\le \BP(\eta (B)=0)P(\eta (B')=0).$$
Therefore for this process we obtain $\BP(\eta_1(B)=0)\le \exp(-\BE\eta(B)).$
\ \ $\square$
\end{example}
For mixed sampled point processes on general spaces we get the following comparison result.
\begin{proposition}\label{mixed-poisson}
Suppose $\eta_1$ is a  mixed sampled point process on  a complete, separable metric space $\mathbb X$  defined by (\ref{mixed}) for which $\tau\in \mathcal S.$ Then
$$\eta_1<_{dcx}\eta_2$$ where  $\eta_2$ denotes a Poisson point process on $\mathbb X$, with the intensity measure $\BE\eta_1.$
\end{proposition}
\proof From Theorem \ref{finiteNA2} we know that $\eta_1$ is $NA$, and from Proposition \ref{dcx-poi} we shall get the conclusion of the present proposition if we  show that for such processes $\eta_1 (B)<_{cx}Y$, where $Y$ denotes a random variable with Poisson distribution $Po(\BE\eta_1(B))$. From the definition of mixed sampled point processes we know that $\eta_1(B)$ is distributed as a random sum $\sum_{i=1}^\tau U_i$, where $(U_i, i\ge 1)$ is an iid sequence of Bernoulli, i.e. 0-1 valued  variables with  success probability $F(B)$. Since  $\tau \in \mathcal S_n$, from the definition of the log-concave ordering $<_{lc}$ in \cite{whitt1985uniform} it follows that  $\tau <_{lc} Po(\BE(\tau))$, which  implies that $\tau<_{cx} Po(\BE(\tau)$, see Theorem 1 in \cite{whitt1985uniform}. It follows that $\sum_{i=1}^\tau U_i<_{cx}\sum_{i=1}^{Po(\BE(\tau))}U_i,$ where $Po(\BE(\tau))$ denotes a Poisson random variable which is independent of $(U_i, i\ge 1)$ (see e.g. \cite{kulik2005dependence}, Corollary 4.5). From this we get $\eta _1(B)<_{cx} Y,$ for $Y$ described above since $\BE\eta_1(B)=\BE(\tau)\BE(U_i)$ and $\sum_{i=1}^{Po(\BE(\tau))}U_i$ has Poisson distribution.
For arbitrary $\tau\in \mathcal S$ we apply weak approximation by $\tau$'s in $\mathcal S_n, n\ge 1.$
\ $\square$

From Proposition \ref{momentrd} we get the following corollary.
\begin{corollary}
Suppose that $\eta$ is a simple mixed sampled point process on $\R^d$, defined by (\ref{mixed}), for which $\tau\in \mathcal S.$ Then
\begin{enumerate}
\item (moment measures) $\BE (\eta (B_1)\cdots\eta  (B_n)))\le \BE (\eta (B_1))\cdots\BE(\eta (B_n)),$\\ for all disjoint bounded Borel sets $B_1,\ldots, B_n$, $n\ge 1,$
\item (void probabilities) $\BP(\eta (B)=0)\le \exp(-\BE\eta(B)),$ for all bounded Borel sets $B$, $n\ge 1$. 
\end{enumerate}
\end{corollary}
Let $\lambda$ be a Radon measure on a locally compact Polish space $\mathbb X$. Let $K$ be a locally trace-class positive contraction on $L_2(\mathbb X, \lambda)$, and  $\eta_K$ be  the determinantal point process generated by $K$. From  Proposition \ref{dcx-poi} we obtain the following corollary.
\begin{corollary}
 Suppose that $\eta_K$ is the determinantal point process described above. Then $$\eta_K<_{dcx}\eta_2$$ where $\eta_2$ denotes a Poisson point process with intensity measure $\BE\eta_K.$
\end{corollary}
\proof  Fix a bounded Borel set $B$. From  Theorem \ref{determinantal-na} we know that $\eta_K$ is $NA$ so in order to get the conclusion of this corollary it is enough (see Proposition \ref{dcx-poi}) to show that $\eta _K(B)<_{cx} Y$, where $Y$ is a random variable with $Po(\BE\eta_K(B))$ distribution.
From Hough et al. \cite{Hough2006}, Proposition 9, we know that $\eta_K (B)$ is distributed as a sum of independent Bernoulli random variables.  From the definition of the log-concave ordering $<_{lc}$ in \cite{whitt1985uniform} it follows that  $\eta_K(B) <_{lc} Y$, where $Y$ is described above, which in turn  implies that $\eta_K (B)<_{cx}Y$, see Theorem 1 in \cite{whitt1985uniform}.
$\square$

The above corollary  for the case of jointly observable sets and $\mathbb X =\R^d$ was observed by Blaszczyszyn and Yogeshwaran in \cite{blaszczyszyn2014comparison}, Proposition 5.3, using a different argument.
  

\section{Some applications} 

Let $\eta$ be a point process on  a complete, separable metric space $\mathbb X$.  Using Chebyshev's inequality we have
$$\BP(|\eta (B)-\BE(\eta(B))| \ge \epsilon) \le \BV (\eta (B)/(\epsilon ^2),$$ for all bounded Borel sets  $B$ and $\epsilon >0$. 

Similarly, using Chernoff's bound,
$$
\BP (\eta (B)-\BE (\eta (B))\ge \epsilon)\le e^{-t(\BE(\eta(B))+\epsilon)}\BE(e^{t\eta (B)}),
$$
for any $t,a >0$, and the upper bounds can be replaced by the values taken from the dominating in $<_{dcx}$ process. If $\eta$ is determinantal or $NA$ mixed sampled point processes, the corresponding Poisson processes can be used to obtain  upper bounds and  concentration inequalities using Proposition \ref{mixed-poisson}.

Similarly, we have for all bounded Borel sets  $B$ and $\epsilon , t>0$. 
$$\BP(\BE(\eta(B))-\eta(B)\ge \epsilon)\le e^{t(\BE(\eta(B))-\epsilon)}\BE (e^{-t\eta(B)}).$$

Using Corollary 2 from \cite{Christofides2004} we can get, from negative association of $\eta$, Kolmogorov type inequalities.
\begin{corollary}
Suppose that $\eta$ is a mixed sampled point process  on  a complete, separable metric space $\mathbb X$, for which $\tau\in \mathcal S$.  Then for any increasing sequence $b_k, k\ge 1$ of positive numbers, any collection of disjoint bounded Borel sets $B_1,\ldots, B_n\in {\mathcal X}$, and $\epsilon >0$
\begin{enumerate}
\item
$$
\BP(\max_{k\le n}|\frac{1}{b_k}\sum_{i=1}^k(\eta (B_i)-\BE(\eta (B_i))|\ge \epsilon)\le 
8\epsilon ^{-2}\sum_{i=1}^n\frac{\BV(\eta (B_i))}{b_i^2},
$$
\item for any integer $m<n$
$$
\BP(\max_{m\le k\le n}|\frac{1}{b_k}\sum_{i=1}^k(\eta (B_i)-\BE(\eta (B_i))|\ge \epsilon)\le 
32\epsilon ^{-2}(\sum_{i=m+1}^n\frac{\BV(\eta (B_i))}{b_i^2}+\sum_{i=1}^m\frac{\BV(\eta (B_i))}{b_m^2}).
$$
\end{enumerate}
\end{corollary}


\bibliography{proba2}
\bibliographystyle{plain}
\end{document}